\documentclass[12pt, twoside, eqno]{article}
\usepackage{latexsym}
\usepackage{amssymb}
\usepackage{amsfonts}
\begin{document}\bigskip

\noindent{\Large\bf On hypersemigroups}\bigskip

\noindent{\bf Niovi Kehayopulu and Michael Tsingelis}\bigskip

\noindent March 4, 2015\bigskip{\small

\noindent{\bf Abstract.} In this paper we show the way we pass from 
semigroups (without order) to hypersemigroups. Moreover we show that, 
exactly as in semigroups, in the results of hypersemigroups based on 
right (left) ideals, quasi-ideals and bi-ideals, points do not play 
any essential role, but the sets, which shows their pointless 
character. The aim of writing this paper is not just to add a 
publication on hypersemigroups but, mainly, to publish a paper which 
serves as an example to show what an hypersemigroup is and give the 
right information concerning this structure. \medskip

\noindent{\bf AMS 2010 Subject Classification:} 20M99\medskip

\noindent{\bf Keywords:} Hypersemigroup; regular; left (right) ideal; 
intra-regular; bi-ideal; quasi-ideal }

\section {Introduction and prerequisites} This paper serves as an 
example to show how similar is the theory of semigroups (without 
order) with the theory of hypersemigroups (called semihypergroups as 
well). A semigroup $(S,.)$ is called {\it regular} if for every $a\in 
S$ there exists an element $x\in S$ such that $a=axa$. This is 
equivalent to saying that $a\in aSa$ for every $a\in S$ or 
$A\subseteq ASA$ for every $A\subseteq S$. A semigroup $(S,.)$ is 
called {\it intra-regular} [1] if for every $a\in S$ there exist 
$x,y\in S$ such that $a=xa^2y$ that is if $a\in Sa^2S$ for every 
$a\in S$ or $A\subseteq SA^2S$ for every $A\subseteq S$. A nonempty 
subset $A$ of $S$ is called a {\it right} (resp. {\it left}) {\it 
ideal} of $(S,.)$ if $AS\subseteq A$ (resp. $SA\subseteq A)$. A 
nonempty subset $B$ of $S$ is called a {\it bi-ideal} of $S$ if 
$BSB\subseteq B$.
The concept of regular rings was introduced by J.v. Neumann [8]. L. 
Kov\'acs characterized the regular rings as rings satisfying the 
property $A\cap B=AB$ for every right ideal $A$ and every left ideal 
$B$ of $S$, where $AB$ is the set of all finite sums of the form 
$\sum\limits_{{a_i} \in A,\,\,{b_i} \in B} {{a_i}} {b_i}$ [3]. K. 
Iseki studied the same for semigroups, and showed that a semigroup 
$S$ is regular if and only if for every right ideal $A$ and every 
left ideal $B$ if $S$, we  have $A\cap B=AB$ [2]. A semigroup $S$ is 
intra-regular if and only if for any right ideal $A$ and every left 
ideal $B$ of $S$, we have $A\cap B\subseteq BA$ [6]. In addition, a 
semigroup $S$ is regular, intra-regular and both regular and 
intra-regular if and only if for every right ideal $X$, every left 
ideal $Y$ and every bi-ideal $B$ of $S$ we have $X\cap B\cap 
Y\subseteq XBY$, $Y\cap B\cap X\subseteq YBX$ or $X\cap B\cap 
Y\subseteq BXY$, respectively [4, 5]. We examine these results for 
hypersemigroups in an attempt to show how similar is the theory of 
semigroups with the theory of hypersemigroups. We tried to use sets 
instead of elements to show that, exactly as in semigroups, for the 
results on hypersemigroups based on ideals given in this paper, 
points do not play any essential role, but the sets which shows their 
pointless character. We show that the corresponding results on 
semigroups are also hold as application of the results of this paper. 
But, though the concept of hypersemigroup generalizes the concept of 
a semigroup, if we want to get a result on an hypersemigroup, we 
never work directly on the hypersemigroup. Exactly as in 
Gamma-semigroups, we have to examine it first for a semigroup. The 
present paper serves as an example to justify what we say.

For an hypersemigroup $H$ we have two operations. One of them is the 
operation between the elements of $H$ and the other between the 
nonempty subsets of $H$. If we show these two operations with the 
same symbol (as we have often seen in the bibliography), a great 
confusion erases. So we have to show these two operations by 
different symbols. As the operation between the elements of $H$ is 
(mostly) denoted by ``$o$", we denote the operation between the 
nonempty subsets of $H$ by the symbol ``$*$". For convenience, let us 
give the following definitions-notations: Let $H$ be a nonempty set 
and ${\cal P}^*(H)$ the set of all nonempty subsets of $H$. Any 
mapping$$\circ : H\times H \rightarrow {\cal P}^*(H) \mid (a,b) 
\rightarrow a\circ b$$is called an {\it hyperoperation} on $H$ and 
the pair $(H,\circ)$ is called and {\it hypergroupoid}.\\Let ``$*$" 
be the operation on ${\cal P}^*(H)$ (induced by the operation of $H$) 
defined by:$$* : {\cal P}^*(H)\times {\cal P}^*(H) \rightarrow {\cal 
P}^*(H) \mid (A,B) \rightarrow A*B:=\bigcup\limits_{(a,b) 
\in\,A\times B} {(a\circ b)}.$$We can easily show that the operation 
``$*$" is well defined.\medskip

\noindent{\bf Proposition 1.} {\it If $(H,\circ)$ is an hypergroupoid 
then, for any $A,B,C,D\in {\cal P}^*(H)$, we have

$(1)$ $A\subseteq B
\mbox { and } C\subseteq D\;\Longrightarrow\; A*C\subseteq B*D$ and

$(2)$ $A\subseteq B \Longrightarrow A*C\subseteq B*C \mbox { and } 
C*A\subseteq C*B.$}\medskip

\noindent{\bf Definition 2.} An hypergroupoid $(H,\circ)$ is called 
{\it hypersemigroup} if $$\{x\}*(y\circ z)=(x\circ y)*\{z\}$$for all 
$x,y,z\in H$.\medskip

Since $\{x\}*(y\circ z)=\bigcup\limits_{\scriptstyle a \in \{ x\} 
\hfill\atop
\scriptstyle b \in (y\circ z)\hfill} {(a\circ b) = \bigcup\limits_{b 
\in (y\circ z)} {(x\circ b)} }$ and

$(x\circ y)*\{z\}=\bigcup\limits_{\scriptstyle a\in \{ x\circ y\} 
\hfill\atop
\scriptstyle b \in \{z\}\hfill} {(a\circ b) = \bigcup\limits_{a\in 
(x\circ y)} {(a\circ z)} },$ an hypergroupoid $H$ is an 
hypersemigroup if and only if $\bigcup\limits_{b \in (y\circ z)} 
{(x\circ b)}=\bigcup\limits_{a\in (x\circ y)} {(a\circ z)}$ for for 
every $x,y,z\in H$.\medskip

In addition, we denote by ${\cal P}^*(H)$ the set of nonempty subsets 
of $H$ and, for an hypersemigroup, we certainly have to show that the 
operation ``$*$" on ${\cal P}^*(H)$ satisfies the associativity 
relation which leads to the associativity relation on any finite 
sequence $A_1,A_2, ... , A_n$ of elements of ${\cal P}^*(H)$. Unless 
this associativity relation, in an expression of the form $A_1*A_2 
***\, A_n$ of elements of ${\cal P}^*(H)$ it is not known where to 
put the parenthesis, and so all the results we get are without 
sense.

We use the terms left (right) ideal, bi-ideal, quasi-ideal instead of 
left (right) hyperideal, bi-hyperideal, quasi-hyperideal and so on, 
and this is because in this structure there are not two kind of left 
ideals, for example, to distinguish them as left ideal and left 
hyperideal. The left ideal in this structure is that one which 
corresponds to the left ideal of semigroups.
\section{Main results}
\noindent{\bf Remark 3.} If $H$ is an hypergroupoid then, for every 
$x,y\in H$, we have$$\{x\}*\{y\}=x\circ y.$$Indeed, 
$\{x\}*\{y\}=\bigcup\limits_{\scriptstyle u \in \{ x\} \hfill\atop
\scriptstyle v \in \{ y\} \hfill} {(u\circ v)}=x\circ y.$\\
Let us give first an example of an hypersemigroup.\medskip

\noindent{\bf Example 4.} If $(S,.)$ is a semigroup and ``$\circ$" 
the mapping of $S\times S$ into ${\cal P}^*(S)$ defined by$$\circ : 
S\times S \rightarrow {\cal P}^*(S) \mid (x,y) \rightarrow x\circ 
y:=\{xy\},$$then $(S,\circ)$ is an hypersemigroup. 
Indeed,\begin{eqnarray*}\{x\}*(y\circ z)&=&\{x\}*\{yz\} \mbox { (the 
operation ``$*$" is well defined) }\\&=&x\circ (yz) \mbox { (by 
Remark 3)}\\&=&\{x(yz)\}=\{(xy)z\}=(xy)\circ z\\&=&\{xy\}*\{z\}\mbox 
{ (by Remark 3)}\\&=&(x\circ y)*\{z\}\mbox { (the operation ``$*$" is 
well defined)}.\end{eqnarray*}The following proposition, though 
clear, plays an essential role in the theory of hypersemigroups.
\medskip

\noindent{\bf Proposition 5.} {\it Let $(H,\circ)$ be an 
hypergroupoid, $x\in H$ and $A,B\in {\cal P}^*(H)$. Then we have the 
following:

$(1)$ $x\in A*B$ $\Longleftrightarrow$ $x\in a\circ b$ for some $a\in 
A$, $b\in B$.

$(2)$ If $a\in A$ and $b\in B$, then $a\circ b\subseteq 
A*B$.}\medskip

\noindent{\bf Proposition 6.} {\it Let $(H,\circ)$ be an 
hypergroupoid and $A_i,B\in {\cal P}^*(H)$, $i\in I$. Then we have 
the following:

$(1)$ $(\bigcup\limits_{i \in I} {{A_i}} )*B = \bigcup\limits_{i \in 
I} {({A_i}} *B)$.\smallskip

$(2)$ $B*(\bigcup\limits_{i \in I} {{A_i}} ) = \bigcup\limits_{i \in 
I} {(B*{A_i}} )$.}\medskip

\noindent{\bf Proof.} (1) Let $x\in (\bigcup\limits_{i \in I} {{A_i}} 
)*B$. By Proposition 5(1), we have $x\in a\circ b$ for some $a\in 
(\bigcup\limits_{i \in I} {{A_i}} )$, $b\in B$. Since $a\in A_j$ for 
some $j\in I$ and $b\in B$, by Proposition 5(2), we have $a\circ 
b\subseteq A_j*B\subseteq \bigcup\limits_{i \in I} {({A_i}} *B)$.
Let now $x\in \bigcup\limits_{i \in I} {({A_i}} *B)$. Then $x\in 
A_j*B$ for some $j\in I$. By Proposition 5(1), we have $x\in a\circ 
b$ for some $a\in A_j$, $b\in B$. Since $a\in A_j\subseteq 
\bigcup\limits_{i \in I} {{A_i}} $ and $b\in B$, by Prop. 5(2), we 
have $a\circ b\subseteq (\bigcup\limits_{i \in I} {{A_i}} )*B$. Then 
we get $x\in (\bigcup\limits_{i \in I} {{A_i}} )*B$.\\The proof of 
property (2) is similar.$\hfill\Box$\medskip

\noindent{\bf Proposition 7.} {\it Let $(H,\circ)$ be an 
hypergroupoid, $A_i,B\in {\cal P}^*(H)$, $i\in I$ and 
$\bigcap\limits_{i \in I} {{A_i}}\not=\emptyset$. Then we have the 
following:

$(1)$ $(\bigcap\limits_{i \in I} {{A_i}} )*B\subseteq 
\bigcap\limits_{i \in I} {({A_i}} *B)$.\smallskip

$(2)$ $B*(\bigcap\limits_{i \in I} {{A_i}} )\subseteq 
\bigcap\limits_{i \in I} {(B*{A_i}} )$.}\medskip

\noindent{\bf Proof.} Let us prove the property (1). The proof of (2) 
is similar.\\(1) First of all, $\bigcap\limits_{i \in I} {{A_i}}\in 
{\cal P}^*(H)$ (since $\bigcap\limits_{i \in I} 
{{A_i}}\not=\emptyset$). Let $x\in (\bigcap\limits_{i \in I} {{A_i}} 
)*B$. Then, by Prop. 5(1), we have $x\in a\circ b$ for some $a\in 
\bigcap\limits_{i \in I} {{A_i}}$, $b\in B$. Since $a\in A_i$ 
$\forall$ $i\in I$ and $b\in B$, by Prop. 5(2), we have $a\circ 
b\subseteq A_i*B$ $\forall$ $i\in I$. Thus we get $x\in a\circ 
b\subseteq \bigcap\limits_{i \in I} {({A_i}} 
*B)$.$\hfill\Box$\medskip

\noindent{\bf Proposition 8.} {\it Let $(H,\circ)$ be an 
hypersemigroup and $A,B,C\in {\cal P}^*(H)$. Then we have 
$$(A*B)*C=A*(B*C)=\bigcup\limits_{(a,b,c)\in A\times B\times C} 
{{\Big(}(a\circ b)*\{ c\}\Big) }.$$}{\bf Proof.} In 
fact,\begin{eqnarray*}(A*B)*C&=&{\Big(}\bigcup\limits_{(a,b) 
\in\,A\times B} {(a\circ b)}{\Big)}*C \mbox { (the operation ``$*$" 
is well defined)}\\&=&\bigcup\limits_{(a,b) \in\,A\times B} 
{\Big(}{(a\circ b)}*C{\Big)} \mbox { (by Proposition 
6(1))}\\&=&\bigcup\limits_{(a,b) \in\,A\times B} {\Big(}{(a\circ 
b)}*\bigcup\limits_{c \in C} {\{ c\} }{\Big)}.\end{eqnarray*}By 
Proposition 6(2), we have$$(a\circ b)*\bigcup\limits_{c \in C} {\{ 
c\} }=\bigcup\limits_{c \in C} {{\Big(}(a\circ b)*\{ c\}\Big) 
}.$$Thus we have$$(A*B)*C=\bigcup\limits_{(a,b) \in\,A\times B} 
\;\;{\Bigg(}\bigcup\limits_{c \in C} {{\Big(}(a\circ b)*\{ 
c\}\Big)}{\Bigg)}=\bigcup\limits_{(a,b,c)\in A\times B\times C} 
{{\Big(}(a\circ b)*\{ c\}\Big) }.$$On the other hand,
\begin{eqnarray*}A*(B*C)&=&A*{\Big(}\bigcup\limits_{(b,c) 
\in\,B\times C} {(b\circ c)}{\Big)}=\bigcup\limits_{(b,c) 
\in\,B\times C}{\Big(}A*{(b\circ c)}{\Big)}\\&=&\bigcup\limits_{(b,c) 
\in B\times C} {{\Bigg(}{\Big(}\bigcup\limits_{a \in A} {\{ a\} 
{\Big)}*(b\circ c){\Bigg)}} }\\&=&\bigcup\limits_{(b,c) \in B\times 
C} {{\Bigg(}\bigcup\limits_{a \in A}{\Big(} {\{ a\} *(b\circ 
c){\Big)}} }{\Bigg)}\\&=&\bigcup\limits_{(a,b,c) \in A\times B\times 
C}{\Big(}\{a\}*(b\circ c){\Big)}.\end{eqnarray*}Since $(H,*)$ is a 
hypersemigroup, we have $(a\circ b)*\{c\}=\{a\}*(b\circ c)$, thus we 
have $(A*B)*C=A*(B*C).$ Moreover, since $b\circ c=\{b\}*\{c\}$, we 
have $\{a\}*(b\circ c)=\{a\}*\{b\}*\{c\}$. $\hfill\Box$\medskip

\noindent According to Proposition 8, we 
write$$(A*B)*C=A*(B*C)=A*B*C.$$For convenience, it is no harm to 
prove the following proposition as well, needed in the definition of 
intra-regular hypersemigroups.\medskip

\noindent {\bf Proposition 9.} {\it Let $(H,\circ)$ be an 
hypersemigroup and $A,B,C,D\in {\cal P}^*(H)$. Then we have 
\begin{eqnarray*}(A*B*C)*D&=&A*(B*C*D)=(A*B)*(C*D)\\&=&\bigcup\limits_{(a,b,c,d)\in 
A\times B\times C\times D} {\Bigg(}{{\Big(}(a\circ b)*\{ 
c\}{\Big)}*\{d\}{\Bigg)} }. \end{eqnarray*}}{\bf Proof.} By 
Proposition 8, we have 
$$(A*B*C)*D={\Big(}(A*(B*C){\Big)}*D=A*{\Big(}(B*C)*D{\Big)}=A*(B*C*D)$$
and$$(A*B*C)*D={\Big(}(A*B)*C{\Big)}*D=(A*B)*(C*D).$$So we 
write$$(A*B*C)*D=A*(B*C*D)=(A*B)*(C*D)=A*B*C*D.$$Let us give the form 
of the elements of $A*B*C*D$. We have
\begin{eqnarray*}(A*B*C)*D&=&{\Bigg(}\bigcup\limits_{(a,b,c)\in 
A\times B\times C} {{\Big(}(a\circ b)*\{ c\}\Big){\Bigg)} }*D\mbox { 
(by Prop. 8)}\\&=&\bigcup\limits_{(a,b,c) \in A\times B\times C)} 
{{\Bigg(}{\Big(}(a\circ b)*\{ c\}{\Big)}*D{\Bigg)}} \mbox { (by Prop. 
6(1))}\\&=&\bigcup\limits_{(a,b,c) \in A\times B\times C)} 
{{\Bigg(}{\Big(}(a\circ b)*\{ c\}{\Big)}*\bigcup\limits_{d \in D} {\{ 
d\} }{\Bigg)}}.\end{eqnarray*}By Proposition 6(2), we have

$${\Big(}(a\circ b)*\{ c\}{\Big)}*\bigcup\limits_{d \in D} {\{ d\} 
}=\bigcup\limits_{d \in D} {{\Bigg(}{\Big(}(a\circ b)*\{ c\} 
{\Big)}*\{ d\}{\Bigg)}}.$$Thus we 
have\begin{eqnarray*}(A*B*C)*D&=&\bigcup\limits_{(a,b,c) \in A \times 
B \times C}{\Bigg(} {\bigcup\limits_{d \in D} {{{\Bigg(}{\Big(}(a 
\circ b)*\{ c\}{\Big)}*\{ d\}{\Bigg)}{\Bigg)}}} 
}\\&=&\bigcup\limits_{(a,b,c,d) \in A \times B \times C \times D} 
{{\Bigg(}{\Big(}(a \circ b)*\{ c\} {\Big)}*\{ d\}{\Bigg)}}. 
\end{eqnarray*}
On the other hand,\begin{eqnarray*}{\Big(}(a\circ 
b)*\{c\}{\Big)}*\{d\}&=&{\Big(}\{a\}*(b\circ 
c){\Big)}*\{d\}=\{a\}*{\Big(}(b\circ 
c)*\{d\}{\Big)}\\&=&\{a\}*\{b\}*\{c\}*\{d\}.\end{eqnarray*} 
$\hfill\Box$

By Proposition 8, using induction (exactly as in Proposition 9), for 
any finite family $A_1,A_2, ..., A_n$ of elements of ${\cal P}^*(H)$, 
we have$$A_1\times A_2\times ..... \times 
A_n=\bigcup\limits_{({a_1},{a_2}...{a_n} )\in {A_1} \times {A_2} 
\times ... \times {A_n}}{\Big(} {\{ {a_1}\}*\{ {a_2}\}* ,... * \{ 
{a_n}} \}{\Big)}.$$ $\hfill\Box$

We are ready now to give the correct definition of regular and 
intra-regular hypersemigroups and give the characterizations of them 
which correspond to the characterization of regular and intra-regular 
semigroups mentioned in the introduction.\medskip

\noindent{\bf Definition 10.} An hypersemigroup $(H,\circ)$ is called 
{\it regular} if for every $a\in H$ there exists $x\in H$ such that 
$$a\in (a\circ x)*\{a\}.$$
\noindent{\bf Proposition 11.} {\it Let $(H,\circ)$ be an 
hypersemigroup. The following are equivalent:

$(1)$ H is regular.

$(2)$ $a\in \{a\}*H*\{a\}$ for every $a\in H$.

$(3)$ $A\subseteq A*H*A$ for every $A\in {\cal P}^*(H)$.}\medskip

\noindent{\bf Proof.} $(1)\Longrightarrow (2)$. Let $a\in H$. Since 
$a$ is regular, there exists $x\in H$ such that $a\in (a\circ 
x)*\{a\}$. Since $a\in \{a\}$ and $x\in H$, by Proposition 5(2), we 
have $a\circ x\subseteq \{a\}*H$. Since $a\circ x\subseteq \{a\}*H$ 
and $\{a\}\subseteq \{a\}$, by Proposition 1, we have $(a\circ 
x)*\{a\}\subseteq {\Big(}\{a\}*H{\Big)}*\{a\}$. Thus we have $a\in 
\{a\}*H*\{a\}$.\\$(2)\Longrightarrow (3)$. Let $A\in {\cal P}^*(H)$ 
and $a\in A$. By (2), we have $a\in \{a\}*H*\{a\}$. Since  
$\{a\}\subseteq A$, by Proposition 1, we have $\{a\}*H\subseteq A*H$. 
Since $\{a\}*H\subseteq A*H$ and $\{a\}\subseteq A$, by Proposition 
1, we have ${\Big(}\{a\}*H{\Big)}*\{a\}\subseteq (A*H)*A.$ Then we 
have $$a\in \{a\}*H*\{a\}={\Big(}\{a\}*H{\Big)}*\{a\}\subseteq 
(A*H)*A=A*H*A$$ and $a\in A*H*A$.\\$(3)\Longrightarrow (1)$. Let 
$a\in H$. By (3), we have $\{a\}\subseteq \{a\}*H*\{a\}$. By 
Proposition 8, $\{a\}*H*\{a\}=\bigcup\limits_{h \in H} 
{{\Big(}(a\circ h)} *\{ a\} {\Big)}$. Then there exists $x\in H$ such 
that $a\in (a\circ x)*\{a\}$, and $H$ is 
regular.$\hfill\Box$\medskip

\noindent{\bf Proposition 12.} {\it If H is an hypergroupoid, then we 
have

$(1)$ $H*A\subseteq H$ and $A*H\subseteq H$ for any $A\in{\cal 
P}^*(H)$.

$(2)$ $H*H\subseteq H$.}\medskip

\noindent{\bf Proof.} (1) Let $A\in{\cal P}^*(H)$ and
$t\in H*A$. Then, by Proposition 5(1), $t\in u\circ v$ for some $u\in 
H$, $v\in A$. Since $u,v\in H$, we have $u\circ v\in {\cal P}^*(H)$, 
that is $u\circ v\subseteq H$, thus we have $t\in H$. Similarly, 
$A*H\subseteq H$. The property (2) is an immediate consequence of 
(1). $\hfill\Box$\medskip

\noindent{\bf Definition 13.} Let $(H,\circ)$ be an hypergroupoid. A 
nonempty subset $A$ of $H$ is called a {\it left} (resp. {\it right}) 
ideal of $H$ if $H*A\subseteq A$ (resp. $A*H\subseteq A)$. A subset 
of $H$ which is both a left and a right ideal of $H$ is called an 
{\it ideal} (or {\it two-sided ideal}) of $H$.

For every $A\in {\cal P}^*(H)$, the set $A*H$ is a right ideal of 
$H$, $H*A$ is a left ideal of $H$ and $H*A*H$ is an ideal of $H$. So 
$H$ is a right ideal, left ideal and an ideal of $H$.

For any nonempty subset $A$ of $H$, we denote by $R(A)$, $L(A)$ and 
$I(A)$ the right ideal, left ideal and the ideal of $H$, 
respectively, generated by $A$, that is the least with respect to the 
inclusion relation left, right or ideal of $H$ containing $A$. 
Exactly as in semigroups, one can easily prove that, for any $A\in 
{\cal P}^*(H)$, we have

$R(A)=A\cup (A*H)$, $L(A)=A\cup (H*A)$, and

$I(A)=A\cup (H*A)\cup (A*H)\cup (H*A*H)$. \medskip

\noindent{\bf Proposition 14.} {\it Let $(H,\circ)$ be an 
hypergroupoid. If $A$ is a left (resp. right) ideal of H, then for 
every $h\in H$ and every $a\in A$, we have $h\circ a\subseteq A$ 
(resp. $a\circ h\subseteq A$). ``Conversely", if $A$ is a nonempty 
subset of $H$ such that $h\circ a\subseteq A$ (resp. $a\circ 
h\subseteq A$) for every $h\in H$ and every $a\in A$, then the set 
$A$ is a left (resp. right) ideal of H}.\medskip

\noindent{\bf Proof.} $\Longrightarrow$. Let $A$ be a a left ideal of 
$H$, $h\in H$ and $a\in A$. Since $h\in H$ and $a\in A$, by 
Proposition 5(2), we have $h\circ a\subseteq H*A$. Since $A$ is a 
left ideal of $H$, we have $H*A\subseteq A$. Thus we have $h\circ 
a\subseteq A$.\\$\Longleftarrow$. Let $A$ be a nonempty subset of $H$ 
such that $h\circ a\subseteq A$ for every $h\in H$ and every $a\in 
A$. Then $H*A:=\bigcup\limits_{u \in H,\,\,v \in A} {(u\circ v)} 
\subseteq A$, so $A$ is a left ideal of $H$. $\hfill\Box$\medskip

\noindent{\bf Proposition 15.} {\it Let H be an hypersemigroup. If A 
is a right ideal and B a left ideal of H, then $A\cap 
B\not=\emptyset$.}\medskip

\noindent{\bf Proof.} Take an element $a\in A$ and an element $b\in 
B$ $(A,B\not=\emptyset)$. Since $\{a\}\subseteq A$ and 
$\{b\}\subseteq B$, by Proposition 1, we have $\{a\}*\{b\}\subseteq 
A*B$. By Remark 3, $\{a\}*\{b\}=a\circ b$. By Proposition 5(2) and 
Proposition 1, we have$$a\circ b\subseteq A*B\subseteq A*H\subseteq A 
\mbox { and } a\circ b\subseteq A*B\subseteq H*B\subseteq B,$$so 
$a\circ b\subseteq A\cap B$. Since $a\circ b\in {\cal P}^*(H)$, we 
have $(a\circ b)\not=\emptyset$, then $A\cap 
B\not=\emptyset$.\medskip

\noindent{\bf Theorem 16.} (cf. also [7]) {\it An hypersemigroup H is 
regular if and only if for every right ideal A and every left ideal B 
of H, we have $$A\cap B=A*B \mbox { (equivalently},  A\cap B\subseteq 
A*B).$$}{\bf Proof.} $\Longrightarrow$. Let $A$ be a right ideal and 
$B$ a left ideal of $H$. Since $H$ is regular, for the element $A\cap 
B$ of ${\cal P}^*(H)$, we have \begin{eqnarray*}
A\cap B&\subseteq&{\Big(}(A\cap B)*H{\Big)}*(A\cap B) \mbox { (by 
Prop. 11(1)}\Rightarrow \mbox {(3) and Prop. 8)}\\&\subseteq& 
(A*H)*(A\cap B)\subseteq (A*H)*B \mbox { (by Prop. 1)}\\&=&A*(H*B) 
\mbox { (by Prop. 8) }\\&\subseteq&A*B\subseteq (A*H)\cap (H*B) \mbox 
{ (by Prop. 1)}\\&\subseteq&A\cap B.\end{eqnarray*}Thus we have 
$A\cap B=A*B$.\\
$\Longleftarrow$. Let $A\subseteq H$. By hypothesis and Propositions 
6, 8, 9, 12 and 1, we have\begin{eqnarray*}A&\subseteq& R(A)\cap 
L(A)=R(A)*L(A)={\Big(}A\cup (A*H){\Big)}*{\Big(}A\cup (H*A){\Big)}\\
&=&(A*A)\cup {\Big(}(A*H)*A {\Big)}\cup {\Big(}A*(H*A) {\Big)}\cup 
{\Big(}(A*H)*(H*A){\Big)}\\&=&(A*A)\cup (A*H*A)\cup 
(A*H*H*A)\\&=&(A*A)\cup (A*H*A).\end{eqnarray*}Then, by the same 
propositions, we have\begin{eqnarray*}A*A&=&{\Big(}(A*A)\cup 
(A*H*A){\Big)}*A=(A*A*A)\cup (A*H*A*A)\\&\subseteq& 
A*H*A.\end{eqnarray*}Hence we obtain $A\subseteq A*H*A$ and, by 
Proposition $11(3)\Rightarrow (1)$, $H$ is regular. $\hfill\Box$

Using elements instead of sets, a nice proof of the ``$\Leftarrow$" 
part of the above theorem is as follows:\\
$\Longleftarrow$. Let $a\in H$. By hypothesis, we 
have\begin{eqnarray*}a\in 
R(a)*L(a)&=&{\Bigg(}\{a\}\cup{\Big(}\{a\}*H{\Big)}{\Bigg)}*
{\Bigg(}\{a\}\cup{\Big(}H*\{a\}{\Big)}{\Bigg)}\\&=&{\Big(}\{a\}*\{a\}{\Big)}\cup 
{\Big(}\{a\}*H*\{a\}{\Big)}\cup {\Big(}\{a\}*H*H*\{a\}{\Big)}\\
&=&{\Big(}\{a\}*\{a\}{\Big)}\cup 
{\Big(}\{a\}*H*\{a\}{\Big)}.\end{eqnarray*}If
$a\in \{a\}*H*\{a\}$ then clearly, by Proposition $1(2)\Rightarrow 
(1)$, $H$ is regular.\\If $a\in \{a\}*\{a\}$, then $\{a\}\subseteq 
\{a\}*\{a\}$ then, by Proposition 1, $$a\in \{a\}*\{a\}\subseteq 
\{a\}*\{a\}*\{a\}\subseteq \{a\}*H*\{a\},$$ so $H$ is regular. 
$\hfill\Box$\medskip

\noindent{\bf Corollary 17.} {\it A semigroup $(S,.)$ is regular if 
and only if for every right ideal A and every left ideal B of S, we 
have $$A\cap B=AB \mbox { (equivalently, } A\cap B\subseteq 
AB).$$}{\bf Proof.} For the semigroup $(S,.)$, we consider the 
hypersemigroup $(S,\circ)$ given in the Example 4, where the 
operation ``$\circ$" is defined by $x\circ y:=\{xy\}$. For this 
hypersemigroup, we have $A*B=AB$ for every $A,B\in {\cal P}^*(S)$. in 
fact: Let $t\in A*B$. By Proposition 5(1), $t\in a\circ b$ for some 
$a\in A$, $b\in B$. Since $a\circ b=\{ab\}$, we have $t\in \{ab\}= 
\{a\}\{b\}\subseteq AB$. Let now $t\in AB$. Then $t=ab$ for some 
$a\in A$, $b\in B$. By Proposition 5(2), we have $t\in 
\{t\}=\{ab\}=a\circ b\subseteq A*B.$\\
$\Longrightarrow$. Let $(S,.)$ be regular, $A$ a right ideal and $B$ 
a left ideal of $(S,.)$. The hypersemigroup $(S,\circ)$ is regular. 
Indeed: Let $a\in S$. Since $S$ is regular, there exists $x\in S$ 
such that $a=(ax)a$. Then we 
have\begin{eqnarray*}\{a\}&=&\{(ax)a\}=(ax)\circ a=\{ax\}*\{a\} \mbox 
{ (by Remark 3)}\\&=&(a\circ x)*\{a\} \mbox { (the operation ``$*$" 
is well defined)}.\end{eqnarray*}Then $a\in (a\circ x)*\{a\}$, so 
$(S,\circ)$ is regular. The set $A$ is a right ideal of $(S,\circ)$. 
Indeed: Let $s\in S$, $a\in A$. Since $A$ is a right ideal of 
$(S,.)$, we have $as\in A$, then $a\circ s=\{as\}\subseteq A$. 
Similarly, $B$ is a left ideal of $(S,\circ)$. By Theorem 16, we have 
$A\cap B=A*B$. Since $A*B=AB$, we have $A\cap 
B=AB$.\\$\Longleftarrow$. Suppose $A\cap B\subseteq AB$ for any right 
ideal $A$ and any left ideal $B$ of $(S,.)$. Then for every right 
ideal $A$ and every left ideal $B$ of $(S,\circ)$, we have $A\cap 
B\subseteq A*B$. Indeed: If $A$ is a right ideal of $(S,\circ)$, 
$a\in A$ and $s\in S$, then $as\in \{as\}=a\circ s\subseteq A$, so 
$A$ is a right ideal of $(S,.)$. If $B$ is a left ideal of 
$(S,\circ)$, then it is a left ideal of $(S,.)$ as well. By 
hypothesis, we have $A\cap B\subseteq AB$. Since $AB=A*B$, we have 
$A\cap B\subseteq A*B$. By Theorem 16, $(S,\circ)$ is regular. Then 
$(S,.)$ is regular. Indeed: Let $a\in S$. Since $(S,\circ)$ is 
regular, there exists $x\in S$ such that $a\in (a\circ x)*\{a\}$. By 
the definition of ``$\circ$" and since the operation ``$*$" is well 
defined, we have\begin{eqnarray*}a\in \{ax\}*\{a\}&=&ax\circ a \mbox 
{ (by Remark 3)}\\&=&\{(ax)a\}.\end{eqnarray*}Then $a=axa$, and 
$(S,.)$ is regular.\medskip

\noindent{\bf Definition 18.} Let $H$ be an hypersemigroup. A 
nonempty subset $B$ of $H$ is called a  {\it bi-ideal} of $H$ if 
$$B*H*B\subseteq B.$$Equivalent Definition: 
$\{a\}*\{h\}*\{b\}\subseteq B$ for every $a,b\in B$ and every $h\in 
H$. In fact: If $B*H*B\subseteq B$, $a,b\in B$ and $h\in H$ then, by 
Propositions 1 and 8, we have $\{a\}*\{h\}*\{b\}\subseteq 
B*H*B\subseteq B$. If $\{a\}*\{h\}*\{b\}\subseteq B$ for every 
$a,b\in B$ and every $h\in H$ then, by Proposition 
8,$$B*H*B=\bigcup\limits_{(x,y,z) \in B\times H\times B} {{\Big(}\{ 
x\} *\{ y\} *\{ z\}}{\Big)}\subseteq B,$$and $B*H*B\subseteq 
B$.\medskip

\noindent{\bf Proposition 19.} {\it An hypersemigroup $H$ is regular 
if and only if for every right ideal $X$, every left ideal $Y$ and 
every bi-ideal $B$ of $H$, we have} $$X\cap B\cap Y\subseteq 
X*B*Y.$${\bf Proof.} $\Longrightarrow$. Let $X$ be a right ideal, $Y$ 
a left ideal and $B$ a  bi-ideal of $H$. By hypothesis and 
Proposition $11(1)\Rightarrow (3)$, we have\begin{eqnarray*}&&X\cap 
B\cap Y\subseteq (X\cap B\cap Y)*H*(X\cap B\cap Y)\\&\subseteq& 
{\Big(}(X\cap B\cap Y)*H*(X\cap B\cap Y){\Big)}*H*{\Big(}(X\cap B\cap 
Y)*H*(X\cap B\cap Y){\Big)}\\&=&{\Big(}(X\cap B\cap 
Y)*H{\Big)}*{\Big(}(H\cap B\cap Y)*H*(X\cap B\cap 
Y){\Big)}*{\Big(}H*(X\cap B\cap Y){\Big)}\\&\subseteq& 
(X*H)*(B*H*B)*(H*Y)\\&\subseteq&X*B*Y.\end{eqnarray*}$\Longleftarrow$. 
Let $X$ be a right ideal and $Y$ a left ideal of $H$. By Proposition 
12, $H$ is a bi-ideal of $H$. By hypothesis, we have $$X\cap Y=X\cap 
H\cap Y\subseteq (X*H)*Y\subseteq X*Y.$$
Then, by Theorem 16, $H$ is regular.$\hfill\Box$\medskip

\noindent {\bf Definition 20.} An hypersemigroup $(H,\circ)$ is 
called {\it intra-regular} if, for every $a\in H$, there exist 
$x,y\in H$ such that $$a\in \{x\}*\{a\}*\{a\}*\{y\}.$$On the other 
hand, \begin{eqnarray*}\{x\}*\{a\}*\{a\}*\{y\}&=&(x\circ 
a)*\{a\}*\{y\}= \{x\}*\{a\}*(a\circ y)\\&=&(x\circ a)*(a\circ 
y)=\{x\}*(a\circ a)*\{y\}. \end{eqnarray*}

\noindent {\bf Proposition 21.} {\it Let $(H,\circ)$ be an 
hypersemigroup. The following are equivalent:

$(1)$ H is intra-regular.

$(2)$ $a\in H*\{a\}*\{a\}*H$ for every $a\in H$.

$(3)$ $A\subseteq H*A*A*H$ for every $A\in {\cal P}^*(H)$.}\medskip

\noindent{\bf Proof.} $(1)\Longrightarrow (2)$. Let $a\in H$. Since 
$H$ is intra-regular, there exist $x,y\in H$ such that $a\in 
\{x\}*\{a\}*\{a\}*\{y\}$. Since $\{x\}\subseteq H$, by Proposition 1, 
we have $\{x\}*\{a\}\subseteq H*\{a\}$. Since $\{a\}\subseteq \{a\}$, 
by Proposition 1, we get $$\{x\}*\{a\}*\{a\}\subseteq 
{\Big(}H*\{a\}{\Big)}*\{a\}=H*\{a\}*\{a\}.$$ Since $\{y\}\subseteq 
H$, again by Proposition 1, we have$$\{x\}*\{a\}*\{a\}*\{y\}\subseteq 
{\Big (}H*\{a\}*\{a\}{\Big)}*H=H*\{a\}*\{a\}*H,$$ and $a\in 
H*\{a\}*\{a\}*H$.\smallskip

\noindent$(2)\Longrightarrow (3)$. Let $A\subseteq H$ and $a\in A$. 
Since $\{a\}\subseteq A$, by Proposition 1, we have $H*\{a\}\subseteq 
H*A$, then $H*\{a\}*\{a\}\subseteq (H*A)*A$, and 
$H*\{a\}*\{a\}*H\subseteq 
{\Big(}(H*A)*A{\Big)}*H=H*A*A*H$.\\$(3)\Longrightarrow (1)$. Let 
$a\in H$. By hypothesis and Proposition 9, we have 
\begin{eqnarray*}\{a\}\subseteq 
{\Big(}H*\{a\}*\{a\}{\Big)}*H&=&\bigcup\limits_{(u,v,w,t)\in H\times 
\{a\}\times \{a\}\times H} {\Bigg(}{{\Big(}(u\circ v)*\{ 
w\}{\Big)}*\{t\}{\Bigg)} }\\&=&\bigcup\limits_{(u,t)\in H\times H} 
{\Bigg(}{{\Big(}(u\circ a)*\{ a\}{\Big)}*\{t\}{\Bigg)} 
}.\end{eqnarray*}Then there exist $(x,y)\in H\times H$ such that 
\begin{eqnarray*}a&\in& {\Big(}(x\circ a)*\{a\}{\Big)}*\{y\} 
={\Big(}\{x\}*\{a\}*\{a\}{\Big)}*\{y\}\\&=&\{x\}*\{a\}*\{a\}*\{y\},\end{eqnarray*}so 
$H$ is intra-regular. $\hfill\Box$\medskip

\noindent{\bf Theorem 22.} (see also [7]) {\it Let H be an 
hypersemigroup. Then H is intra-regular if and only if for every 
right ideal A and every left ideal B of H, we have $$A\cap B\subseteq 
B*A.$$}{\bf Proof.} $\Longrightarrow$. Since $H$ is intra-regular, by 
Proposition $21(1)\Rightarrow (3)$, we have$$A\cap B\subseteq 
H*(A\cap B)*(A\cap B)*H\subseteq (H*B)*(A*H)\subseteq 
B*A.$$$\Longleftarrow$. Let $A\subseteq H$. By hypothesis, we 
have\begin{eqnarray*}A&\subseteq& R(A)\cap 
L(A)=L(A)*R(A)={\Big(}A\cup (H*A){\Big)}*{\Big(}A\cup 
(A*H){\Big)}\\&=&(A*A)\cup (H*A*A)\cup (A*A*H)\cup 
(H*A*A*H).\end{eqnarray*}Then we have\begin{eqnarray*}A*A&\subseteq& 
(A*A*A)\cup (H*A*A*A)\cup (A*A*H*A)\cup (H*A*A*H*A)\\&\subseteq& 
(H*A*A*H)\cup (A*A*H),\end{eqnarray*}and $H*A*A\subseteq H*A*A*H$. 
Then we have $A\subseteq (A*A*H)\cup (H*A*A*H),$ from 
which$$A*A\subseteq (A*A*A*H)\cup (A*H*A*A*H)\subseteq H*A*A*H.$$Then 
$A*A*H\subseteq H*A*A*H$, and $A\subseteq H*A*A*H$. By Proposition 
$21(3)\Rightarrow (1)$, $H$ is intra-regular. $\hfill\Box$

Using elements instead of sets, the proof of the ``$\Leftarrow$" part 
of the above theorem is as follows:\\
$\Longleftarrow$. Let $a\in H$. By hypothesis, we 
have\begin{eqnarray*}a\in 
L(a)*R(a)&=&{\Bigg(}\{a\}\cup{\Big(}H*\{a\}{\Big)}{\Bigg)}*
{\Bigg(}\{a\}\cup{\Big(}\{a\}*H{\Big)}{\Bigg)}\\&=&
{\Big(}\{a\}*\{a\}{\Big)}\cup{\Big(} H*\{a\}*\{a\}{\Big)}
\cup{\Big(}\{a\}*\{a\}*H{\Big)}\\&\cup&
{\Big(}H*\{a\}*\{a\}*H{\Big)}.\end{eqnarray*}If $a\in 
{\Big(}H*\{a\}*\{a\}*H{\Big)}$ then, by Proposition 21, $H$ is 
intra-regular. If $a\in \{a\}*\{a\}$ then, by Proposition 1, we have
\begin{eqnarray*}a\in \{a\}&\subseteq&\{a\}*\{a\}\subseteq 
\{a\}*\{a\}*\{a\}\subseteq \{a\}*\{a\}*\{a\}*\{a\}\\&\subseteq& 
H*\{a\}*\{a\}*H.\end{eqnarray*}If $a\in H*\{a\}*\{a\}$, 
then\begin{eqnarray*}a\in \{a\}&\subseteq&H*\{a\}*\{a\}\subseteq 
H*{\Big(}H*\{a\}*\{a\}{\Big)}*\{a\}\\&=&(H*H)*\{a\}*\{a\}*\{a\}
\subseteq H*\{a\}*\{a\}*H.\end{eqnarray*}If $a\in \{a\}*\{a\}*H$, in 
a similar way we get $a\in H*\{a\}*\{a\}*H$, and $H$ is 
intra-regular. $\hfill\Box$\medskip

\noindent{\bf Corollary 23.} A semigroup $(S,.)$ is intra-regular if 
and only if for every right ideal $A$ and every left ideal $B$ of $S$ 
we have $A\cap B\subseteq BA$.\medskip

\noindent{\bf Proof.} As in Corollary 17, we consider the 
hypersemigroup $(S,\circ)$ where $x\circ y=\{xy\}$ for every $x,y\in 
S$. We have already seen that $A$ is a left (resp. right) ideal of 
$(S,.)$ if and only if it is a left (resp. right) ideal of 
$(S,\circ)$ and that $A*B=AB$. It remains to prove that $(S,.)$ is 
intra-regular if and only if $(S,\circ)$ is intra-regular: Let $a\in 
S$. If $(S,.)$ is intra-regular, there exist $x,y\in S$ such that 
$a=xa^2y$. Then we have\begin{eqnarray*}a\in 
\{a\}&=&\{xa^2y\}=\{(xa)(ay)\}=(xa)\circ 
(ay)\\&=&\{xa\}*\{ay\}=(x\circ a)*(a\circ y),\end{eqnarray*} and the 
hypersemigroup $(S,\circ)$ is intra-regular. If $(S,\circ)$ is 
intra-regular, then there exist $x,y\in S$ such 
that\begin{eqnarray*}a\in (x\circ a)*(a\circ 
y)&=&\{xa\}*\{ay\}=(xa)\circ (ay)\\&=&\{xa^2y\},\end{eqnarray*}then 
$a=xa^2y$ and the semigroup $(S,.)$ is intra-regular. 
$\hfill\Box$\medskip

According to Proposition 19, if $(S,.)$ is a semigroup, then for any 
right ideal $X$, any left ideal $Y$ and any bi-ideal of $B$ of $S$, 
we have $X\cap B\cap Y\subseteq X*B*Y.$\\In a similar way as in 
Proposition 19, we can prove the following two propositions which 
also generalize the corresponding results of semigroups mentioned in 
the introduction. \medskip

\noindent{\bf Proposition 24.} {\it An hypersemigroup H is 
intra-regular if and only if for every right ideal X, every left 
ideal Y and every bi-ideal B of H, we have $$X\cap B\cap Y\subseteq 
Y*B*X.$$}\noindent{\bf Proposition 25.} {\it An hypersemigroup $H$  
is both regular and intra-regular if and only if for every right 
ideal X, every left ideal Y and every bi-ideal B of $H$, we have 
$$X\cap B\cap Y\subseteq B*X*Y.$$}{\small
\bigskip

\end{document}